# Polyhedra of Genus 3 with 10 Vertices and Minimal Coordinates

2nd April 2006

### Authors

Stefan Hougardy, Frank H. Lutz, and Mariano Zelke

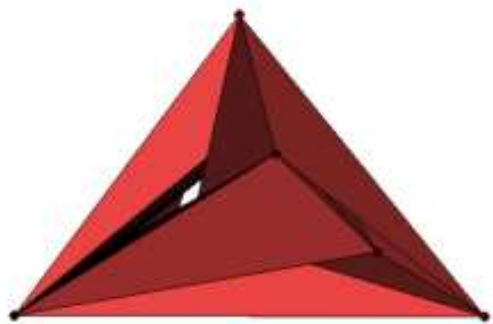

### Description

We give coordinate-minimal geometric realizations in general position for 17 of the 20 vertex-minimal triangulations of the orientable surface of genus 3 in the 5x5x5-cube.

By Heawood's inequality from 1890 [8], every triangulation of a (closed) surface M of Euler characteristic chi(M) has at least

$$n \geq 1/2(7+\sqrt{49-24\cdot chi(M)})$$

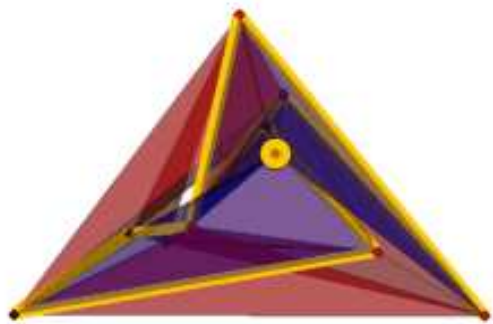

vertices. The tightness of this bound was proved by Jungerman and Ringel [12] for all orientable surfaces (with the exception of the orientable surface of genus 2, where an extra vertex has to be added).

A first vertex-minimal triangulation of the orientable surface of genus 3 with 10 vertices can be found on p. 23 in the book of Ringel on the Map Color Theorem [15]. Polyhedral models for five different 10-vertex triangulations of the orientable surface of genus 3 were given by Brehm [5], [6] and Bokowski and Brehm [2].

Altogether, there are 42426 triangulated surfaces with 10 vertices; cf. [13] and see [14] for a list of facets of the triangulations. In particular, there are exactly 20 combinatorially distinct vertex-minimal 10-vertex triangulations of the orientable surface of genus 3.

By Steinitz' theorem (cf. [17, Ch. 4]), every triangulated 2-sphere is realizable geometrically as the boundary complex of a convex 3-dimensional polytope. For triangulations of orientable surfaces of genus g ≥ 1 it was asked by Grünbaum [7, Ch. 13.2] whether they can always be *realized geometrically as a polyhedron* in $\mathbf{R}^3$, i.e., with straight edges, flat triangles, and without self intersections? In general, the answer turned out to be *NO*: Bokowski and Guedes de Oliveira [3] showed that there is a non-realizable triangulation of the orientable surface of genus 6, and, recently, Schewe [16] was able to extend this result to all surfaces of genus g ≥ 5. However, for surfaces of genus 1 ≤ g ≤ 4 the problem remains open.

Geometric realizations for all 865 vertex-minimal 10-vertex triangulations of the orientable surface of genus 2 were obtained by Bokowski and Lutz [1], [13], based on a random search and geometric intuition.

With a more sophisticated simulated annealing approach, it was also possible to realize surfaces of genus 3:

**Theorem (Hougardy, Lutz, and Zelke [11]):** All 20 vertex-minimal 10-vertex triangulations of the orientable surface of genus 3 can be realized geometrically in $\mathbf{R}^3$.

For most of these examples there even are realizations with rather small coordinates.

**Theorem:** At least 17 of the 20 vertex-minimal 10-vertex triangulations of the orientable surface of genus 3 have realizations in general position in the 5x5x5-cube, but none of the 20 triangulations can be realized in general position in the 4x4x4-cube.

To obtain this result, we completely enumerated for increasing n all sets of 10 vertices in general position in the nxnxn-cube that are compatible with a given triangulation; cf. [9] and [10]. To speed up this enumeration we made use of the symmetry of the nxnxn-cube, enumerated only lexicographic minimal vertex sets, and checked compatibility with a given triangulation for partially generated vertex sets. The search for realizations in the 5x5x5-cube was run (in total) for 2 CPU years on a 3.5 GHz processor. Hereby, roughly 1/5th of the possible vertex sets in the 5x5x5-cube was processed.

Remark: The displayed example Polyhedron_2_10_14542 has one clearly visible hole, while all other tunnels are hidden. In the transparent display of the polyhedron we have highlighted the link of a vertex. The number 14542 indicates the position of the example in the catalog of the 42426 triangulated surfaces with 10 vertices from [14].

**Authors' Addresses**

Stefan Hougardy

>   Humboldt-Universität zu Berlin
>   Institut für Informatik
>   Unter den Linden 6
>   10099 Berlin
>   Germany
>   hougardy@informatik.hu-berlin.de
>   http://www.informatik.hu-berlin.de/~hougardy/

Frank H. Lutz

>   Technische Universität Berlin
>   Fakultät II - Mathematik und Naturwissenschaften
>   Institut für Mathematik, Sekr. MA 3-2
>   Straße des 17. Juni 136
>   10623 Berlin
>   Germany
>   lutz@math.tu-berlin.de
>   http://www.math.tu-berlin.de/~lutz

Mariano Zelke

>   Humboldt-Universität zu Berlin
>   Institut für Informatik
>   Unter den Linden 6
>   10099 Berlin
>   Germany
>   zelke@informatik.hu-berlin.de